\documentclass[a4paper,12pt]{article}
\usepackage[margin=1in]{geometry}
\usepackage[english]{babel} 
\usepackage{epsfig}
\usepackage{fancyhdr} 
\usepackage{amsmath,amssymb, mathtools}
\usepackage{amscd} 
\usepackage[T1]{fontenc} 
\usepackage[utf8]{inputenc} 
\usepackage[usenames,dvipsnames]{xcolor}
\usepackage{graphicx,color,listings}
\usepackage{hologo}
\usepackage{geometry}
\usepackage{rotating}
\usepackage{caption}
\usepackage{amsthm}
\usepackage{stackengine,scalerel}
\usepackage{hyperref}
\hypersetup{citecolor=purple, linkcolor=blue, colorlinks=true}

\def\f{\mathbb F}
\def\fq{\mathbb F_q}

\def\PG{{\rm{PG}}}

\def\pg{{\rm{PG}}}
\DeclareMathOperator\PGO{PGO}

\newcommand{\eps}{\varepsilon}

\theoremstyle{plain}
\usepackage{cleveref}
\newtheorem{thm}{Theorem}[section] 
\newtheorem{cor}[thm]{Corollary} 
\newtheorem{lem}[thm]{Lemma} 
\newtheorem{prop}[thm]{Proposition}

\theoremstyle{definition}
\newtheorem{defn}[thm]{Definition}

\newtheorem{remark}[thm]{Remark}
\newtheorem{construction}[thm]{Construction}

\title{ On $m$-ovoids of $Q^+(7,q)$ with $q$ odd}
\date{}

\author{Sam Adriaensen\thanks{Department of Mathematics and Data Science, Vrije Universiteit Brussel, B-1050 Brussels, Belgium. Email addresses: \{Sam.Adriaensen, Jan.De.Beule, Jonathan.Mannaert\}@vub.be} \and
Jan De Beule\footnotemark[1] \and
Giovanni Giuseppe Grimaldi\thanks{Department of Mathematics and Applications R. Caccioppoli, University of Naples Federico II, Naples, Italy. Email address: giovannigiuseppe.grimaldi@unina.it} \and
Jonathan Mannaert\footnotemark[1]  }

\begin{document}
\maketitle

\begin{abstract}
In this paper, we provide a construction of $(q+1)$-ovoids of the hyperbolic quadric $Q^+(7,q)$, $q$ an odd prime power, 
by glueing $(q+1)/2$-ovoids of the elliptic quadric $Q^-(5,q)$. This is possible by controlling some intersection
properties of (putative) $m$-ovoids of elliptic quadrics. It eventually yields $(q+1)$-ovoids of $Q^+(7,q)$
not coming from a $1$-system. Secondly, for certain values of $q$, we construct line spreads of $\pg(3,q)$ 
that have as many secants to a given elliptic quadric as possible. This is then used to 
construct $m$-ovoids for  $m \in \{ 2,4,6,8,10\}$ in $Q^+(7,3)$. 
\end{abstract}
\thanks{2020 {\em Math. Subj. Class.} 51E20, 05B25}\\
\section{Introduction}

An {\em $m$-ovoid} $\mathcal{O}$ of the elliptic quadric $Q^-(5,q)$ is a set of points intersecting 
any line of $Q^-(5,q)$ in exactly $m$ points.
In 1965, Segre \cite{Segre} proved that if the elliptic quadric $Q^-(5,q)$ has 
a non-trivial $m$-ovoid then $q$ is odd and $m=(q+1)/2$. Secondly, Segre constructed an 
example of a $2$-ovoid of $Q^-(5,3)$. In 2005, Cossidente and Penttila \cite{Coss}
constructed examples of $(q+1)/2$-ovoids of $Q^-(5,q)$ for any $q$ odd, 
generalizing Segre's example.

In this paper, we will construct $m$-ovoids of the hyperbolic quadric $Q^+(7,q)$, $q$ odd, 
as the disjoint union of $m$-ovoids of different elliptic quadrics $Q^-(5,q)$ contained in 
$Q^+(7,q)$. 

Segre's perspective was to study $m$-ovoids of quadrics, i.e. sets of points in a projective
space of which the coordinates satisfy a quadratic equation. Later, the perspective shifted
into using the theory of finite classical polar spaces, a class of geometries comprising quadrics. 
The definition of $m$-ovoids of polar spaces was originally given in \cite{Thas1989}, 
albeit only for polar spaces of rank 2 (also known as \emph{generalized quadrangles}), and
generalized in \cite{ST1994} to polar spaces of arbitrary rank. As stated in \cite{Thas2001}, 
{\em ovoids of polar spaces have many connections with and applications to 
projective planes, circle geometries, generalized polygons, strongly regular 
graphs, partial geometries, semi-partial geometries, codes, and designs}. Ovoids of 
polar spaces are geometrically a generalization of an ovoid in the sense of Tits, defined in 
\cite{Tits1962}, to provide a geometric realization of the Suzuki groups. From the perspective
of algebraic combinatorics, $m$-ovoids are a natural generalization of ovoids, and have strong 
links with strongly regular graphs and two-weight codes see e.g.\ \cite{Bamb,Thas2001}.

The structure of the paper is as follows. Section \ref{Pre} contains preliminaries, including a brief 
introduction to $m$-ovoids and related structures in finite classical polar spaces. 
In Section~\ref{Patt} we determine intersection patterns for
$m$-ovoids of the elliptic quadric $Q^-(2n+1,q)$ with embedded elliptic quadrics $Q^-(2n-1,q)$. For the sake
of completeness, we also present similar results for $m$-ovoids of the symplectic space $W(2n+1,q)$.
In Section~\ref{sec:construction}, we describe how to construct $(q+1)$-ovoids of $Q^+(7,q)$ as a union 
of two disjoint $(q+1)/2$-ovoids of two different embedded elliptic quadrics $Q^-(5,q)$. 
In the final section, we investigate spreads of $\pg(3,q)$ with as many $2$-secants as possible to an elliptic
quadric. This is used eventually to construct $m$-ovoids of $Q^+(7,3)$ for $m \in \{ 2,4,6,8,10 \}$ 
by using 2-ovoids of $Q^-(5,3)$ and the particular line spread. 

\section{Preliminaries}\label{Pre}

\subsection{Finite classical polar spaces}

Polar spaces are incidence geometries originally 
introduced in the seminal work of Tits \cite{Tits} and Veldkamp \cite{Veldkamp}. 
Buekenhout and Shult \cite{BuekenhoutShult} characterized polar spaces 
by the all-or-one axiom. A polar space is a partial linear space with the property that if $P$ is a point and 
$\ell$ a line not incident with $P$, then $P$ is collinear with either all points of 
$\ell$ or with exactly one point of $\ell$. Polar spaces, and many of their substructures,
are quite important and well studied within finite geometry and incidence geometry.
In this paper, we are interested in $m$-ovoids of some particular  
finite classical polar spaces, these are the polar spaces associated to 
sesquilinear and quadratic forms on a vector space. 

Let $\PG(d,q)$ be a $d$-dimensional projective space over the finite field 
$\mathbb{F}_q$ with $q$ elements and suppose that $V=\fq^{d+1}$ is the underlying 
vector space, which is a $(d+1)$-dimensional vector space over $\mathbb{F}_q$.
Let $f$ be a non-degenerate reflexive form, respectively a non-singular quadratic form
on the vector space $\mathbb{F}_q^{d+1}$. We define a {\em classical polar space} $\mathcal{P}$ as 
the incidence structure consisting of the elements of $\pg(d,q)$ with underlying subspaces of $\mathbb{F}_q^{d+1}$
that are totally isotropic, respectively totally singular with respect to the sesquilinear,
respectively quadratic form $f$. The {\em Witt index} of the form $f$ is the dimension of the largest
subspace of $\mathbb{F}_q^{d+1}$ that is totally isotropic, respectively totally singular. So if the
form $f$ has Witt index $r$, the elements of $\mathcal{P}$ are projective subspaces of $\PG(d,q)$
of dimension $0$, $1$, \ldots, $r-1$, hence the polar space $\mathcal{P}$ then consists
of projective points, lines, \ldots, $(r-1)$-dimensional subspaces of $\PG(d,q)$. The elements 
of maximal dimension are called the {\em generators} of $\mathcal{P}$. The geometry $\mathcal{P}$
is naturally embedded in $\PG(d,q)$, and we call $\PG(d,q)$ the {\em ambient projective space}
of $\mathcal{P}$.

Non-degenerate reflexive sesquilinear forms on a vector space $V = \f_q^{d+1}$ are classified, see
e.g. \cite{KL1990}. The three possibilities for such a form $f$ are
\begin{enumerate}
\item $f$ is {\em alternating}, i.e. $f(x,x) = 0$ for all $x \in V$, hence $f(x,y) = - f(y,x)$ for all $x,y \in V$, and $d+1$ is odd, its 
Witt index is $\frac{d+1}{2}$;
\item $f$ is {\em symmetric}, i.e. $f(x,y) = f(y,x)$ for all $x,y \in V$;
\item $f$ is {\em unitary}, i.e. $f(x,y) = f(y,x)^\sigma$, for all $x,y \in V$, with $\sigma$ an involutory automorphism 
of the field $\mathbb{F}_q$ (and hence $q$ is a square).
\end{enumerate}

Symmetric forms in even characteristic do not give rise to a polar space. There are three possibilities for a symmetric 
form $f$ on $\mathbb{F}_q^{d+1}$ in odd characteristic,
\begin{enumerate}
\item $d+1$ is even and $f$ has Witt index $\frac{d+1}{2}$, $f$ is called {\em hyperbolic};
\item $d+1$ is even and $f$ has Witt index $\frac{d-1}{2}$, $f$ is called {\em elliptic};
\item $d+1$ is odd and $f$ has Witt index $\frac{d}{2}$, $f$ is called {\em parabolic}.
\end{enumerate}

There are the same three possibilities for a quadratic form $g$ on $\mathbb{F}_q^{d+1}$ 
according to its Witt index in both even and odd characteristic, and such a form is called 
hyperbolic, respectively, elliptic, parabolic as well. 

Given a quadratic form $g$ on $\mathbb{F}_q^{d+1}$, it is common to call 
the bilinear form $f$ on $\mathbb{F}_q^{d+1}$ defined as $f(x,y) = g(x+y) - g(x) - g(y)$ the {\em polar form} of $g$.
In odd characteristic, $f$ is symmetric, and this relation also defines $g$ uniquely 
from $f$. This is not the case in even characteristic. In that case the polar form $f$ is an
alternating form, which is non-degenerate if and only if $d+1$ is even.

Based on its Witt index, there are two possibilities for a unitary form $f$ on $\mathbb{F}_q^{d+1}$,

\begin{enumerate}
\item $d+1$ is even and $f$ has Witt index $\frac{d+1}{2}$;
\item $d+1$ is odd and $f$ has Witt index $\frac{d}{2}$.
\end{enumerate}

From the synthetic point of view of Tits, Buekenhout and Shult, finite classical polar spaces are
examples of polar spaces. In the synthetic approach, a polar space also has a rank $r$, 
meaning that its subspaces of maximal dimension are indeed $(r-1)$-dimensional projective spaces. 
It then follows from the classification of Veldkamp \cite{Veldkamp} and Tits \cite{Tits} that all 
finite polar spaces of rank at least 3 are classical. The classification of forms on $\fq^{d+1}$
leads to the following list, see also \cite[Theorem 4.3]{Ball}.

\begin{itemize}
\item The {\em elliptic quadric} $Q^-(2n + 1,q)$ of $\mathrm{PG}(2n + 1,q)$ corresponding to the quadratic form
$$
f(x) = X_0X_1 + \ldots + X_{2n-2}X_{2n-1} + g(X_{2n},X_{2n+1}),
$$
where $g$ is a
homogeneous irreducible quadratic form on $\fq^2$.
\item The {\em parabolic quadric} $Q(2n,q)$ in $\mathrm{PG}(2n,q)$ with quadratic form
$$
f(x) = X_0X_1 +\ldots+X_{2n-2}X_{2n-1} +X^2_{2n}.
$$
\item The {\em hyperbolic quadric} $Q^+(2n+1,q)$, in $\mathrm{PG}(2n+1,q)$ with quadratic form
$$
f(x) = X_0X_1 +\ldots+X_{2n}X_{2n+1}.
$$
\item The {\em symplectic polar space} $W(2n+1,q)$ of $\mathrm{PG}(2n+
1,q)$ with alternating form
$$
f(x,y) = X_0Y_1 -X_1Y_0 +\ldots+X_{2n}Y_{2n+1} -X_{2n+1}Y_{2n}.
$$
\item The {\em hermitian variety} $H(n,q^2)$ of $\mathrm{PG}(n,q^2)$ with unitary form 
$$
f(x,y) = X_0^{q+1} + X_1^{q+1} + \ldots + X_n^{q+1}.
$$
\end{itemize}
Recall that the hermitian varieties fall actually apart in two classes depending on 
whether $n$ is odd or even, so we distinguish six different types of finite classical polar
spaces.

From now on, {\em polar space} refers to {\em finite classical polar space}. If $\PG(n,q)$ is the ambient
projective space of a polar $\mathcal{P}$, we say that $\mathcal{P}$ is a polar space in $\PG(n,q)$. 
Let $\mathcal P$ denote a polar space of rank $r$. We call its elements of dimension $k$ simply its \emph{$k$-spaces}.
If $x$ is the number of generators through a fixed $(r-2)$-space of $\mathcal{P}$, we will say that $\mathcal{P}$ has 
{\em parameter $e=\log_q(x-1)$}. \Cref{Table1} lists the different classical polar spaces of rank $r$ and parameter $e$.

\begin{table}[ht]
 \begin{tabular}{|c|c|c|c|c|c|c|}
  \hline
  $\mathcal{P}_{r,e}$ & $Q^+(2r-1,q)$ & $H(2r-1,q)$ & $W(2r-1,q)$ & $Q(2r,q)$ & $H(2r,q)$ & $Q^-(2r+1,q)$\\
  \hline
  $e$ & $0$ & $1/2$ & $1$ & $1$ & $3/2$ & $2$\\
  \hline
 \end{tabular}
 \caption{$\mathcal{P}_{r,e}$ polar space of rank $r$ and parameter $e$}
 \label{Table1}
\end{table}

Let $\mathcal{P}$ be a polar space in $\PG(n,q)$ and $Q$ a point of $\PG(n,q)$. Let $f$ be the underlying
sesquilinear form, or in case $\mathcal{P}$ is associated to a quadratic form $g$, $f$ is the polar form 
of $g$. We say that two points of $\PG(n,q)$ with respective coordinate vectors $x$ and $y$ are orthogonal 
with respect to $\mathcal{P}$ if and only if $f(x,y) = 0$. For a point $Q \in \PG(n,q)$, let $Q^\perp$ denote
all the points of $\PG(n,q)$ orthogonal to $Q$ with respect to $\mathcal{P}$. The set $Q^{\perp}$ is a hyperplane 
of $\PG(n,q)$, and if $Q$ belongs to $\mathcal{P}$ then $Q^{\perp}$ is called the {\em tangent hyperplane at 
$Q$ to $\mathcal{P}$}. For any set $A$ of points, we define $A^{\perp} \coloneqq \cap_{P \in A} P^{\perp}$. 

If $\alpha$ is any subspace of $\PG(n,q)$ then the form $f$ associated to $\mathcal{P}$ will induce a form on $\alpha$ as well
(technically speaking, on the underlying vector space of $\alpha$). This induced form can be degenerate or singular, in which
case it has a non-trivial radical, which induces a subspace $\rho$ of $\PG(n,q)$ contained in $\alpha$. The
intersection of $\mathcal{P}$ with $\alpha$ is then a {\em degenerate polar space}, i.e. a cone with vertex $\rho$ and
base a non-degenerate polar space with associated form $f'$, having the same type as $f$, i.e. alternating, 
symmetric or unitary as well. The following proposition holds, see e.g. \cite{HirsThas}.

\begin{prop}\label{intpolar}
Let $\mathcal{P}_{r,e}$ be a polar space in $\mathrm{PG}(n,q)$.  
\begin{itemize}
\item For any point $P$ of $\mathcal{P}_{r,e}$, the set $P^{\perp} \cap \mathcal{P}_{r,e}$
is a cone with base $\mathcal{P}_{r-1,e}$ and vertex $P$, with $\mathcal{P}_{r-1,e}$ a finite
classical polar space of rank $r-1$ of the same type as $\mathcal{P}_{r,e}$.
\item Let $\ell = \langle P, Q \rangle$ where $P, Q$ are two distinct points of $\mathcal{P}_{r,e}$. Then 
$\ell^{\perp}$ is
an $(n - 2)$-space of $\mathrm{PG}(n, q)$ intersecting $\mathcal{P}_{r,e}$ in a polar space $\mathcal{P}_{r-1,e}$ if
$P \not \in Q^{\perp}$, or in a cone having as vertex the line $\ell$ and as base a
polar space $\mathcal{P}_{r-2,e}$, if $P  \in Q^{\perp}$.
\item Let $\Pi_{k-1}$ be a $(k-1)$-space of $\mathcal{P}_{r,e}$. Then 
$\Pi_{k-1}^\perp$ is
an $(n - k)$-space of $\mathrm{PG}(n, q)$ intersecting $\mathcal{P}_{r,e}$ in a cone having as vertex $\Pi_{k-1}$ and as base a
polar space $\mathcal{P}_{r-k,e}$, where $\mathcal{P}_{r-k,e}$ is of the same type of $\mathcal{P}_{r,e}$.
\end{itemize}
\end{prop}

\begin{prop}[{\cite[Theorem 4.10]{Ball}}]
The number of points of $\mathcal{P}_{r,e}$ is given by
\[
\frac{(q^r-1)(q^{r+e-1}+1)}{(q-1)}.
\]
\end{prop}

The intersection properties of subspaces of the ambient projective space with the polar space play a very important role in the sequel. 
The next two theorems summarize results found in \cite{HirsThas} for quadrics and symplectic polar spaces. We will not consider 
hermitian polar spaces further.

Let $Q^\varepsilon(n,q)$ with $\varepsilon \in \{-1,0,1\}$ denote a quadric which is hyperbolic 
for $\eps = 1$, parabolic for $\eps = 0$, and elliptic for $\eps = -1$.
\begin{thm}[{\cite[Theorem 1.51, Theorem 1.53]{HirsThas}}]
 Suppose that $Q^\eps(n,q)$ is a quadric with $q$ odd. Let $\Sigma$ be any $k$-space of $\pg(n,q)$ 
 and denote $I_1 = Q^\eps(n,q) \cap \Sigma$ and $I_2 = Q^\eps(n,q) \cap \Sigma^\perp$.
 Then both $I_1$ and $I_2$ are quadrics in their respective ambient spaces $\Sigma$ and $\Sigma^\perp$ 
 both with radical $\Pi = \Sigma \cap \Sigma^\perp$. Let $k-t-1$ be the dimension of $\Pi$.
 \begin{itemize}
  \item If $\eps = \pm 1$, then $I_1 = \Pi Q^\delta(t,q)$ for some $\delta \in \{\pm 1, 0\}$, and $I_2 = \Pi Q^{\eps \delta}(n-2k+t-1,q)$.
  \item If $\eps = 0$, then $I_1 = \Pi Q(t,q)$ if and only if $I_2 = \Pi Q^{\pm}(n-2k+t-1,q)$.
 \end{itemize}
\end{thm}

\begin{thm}
Consider the polar space $\mathcal{P} = W(n,q)$, $n$ odd. Let $\Sigma$ be a $k$-space of $\pg(n,q)$. 
Then $\Sigma^\perp$ is an $(n-k-1)$-space. Let $I_1\coloneqq\Sigma \cap \mathcal{P}$ and 
$I_2\coloneqq\Sigma^\perp \cap\mathcal{P}$. Then for some odd $t$ depending on $\Sigma$, 
$I_1=\Pi_{k-t-1}W(t,q)$ and $I_2=\Pi_{k-t-1}W(n - 2k + t - 1, q)$.
\end{thm}

A $k$-space $\Sigma$ of $\PG(n,q)$ will be called {\em hyperbolic}, respectively {\em elliptic, 
parabolic} or {\em symplectic} accordingly to $\Sigma$ meets the polar space 
$\mathcal{P}$ in $Q^+(k,q)$, respectively $Q^-(k,q)$, $Q(k,q)$ or $W(k,q)$.

We refer to \cite{Ball, Cam, HirsThas} for further background on polar spaces.

\subsection{\texorpdfstring{$m$-ovoids and $l$}{{\it m}-ovoids and {\it l}}-systems of finite classical polar spaces}

The central object of this paper are $m$-ovoids.

\begin{defn}
Let $\mathcal{P}_{r,e}$ be a polar space. A set of points $\mathcal{O}$ of $\mathcal{P}_{r,e}$ 
is an {\em $m$-ovoid} if each generator of $\mathcal{P}_{r,e}$ meets $\mathcal{O}$ in $m$ points.
\end{defn}

The following lemma is well known and easily proved, see e.g.\ \cite{Thas1981}.
 
\begin{prop}
Suppose that $\mathcal{O}$ is an $m$-ovoid, then $|\mathcal{O}|=m(q^{r+e-1} + 1)$.
\end{prop}

\begin{lem}[{\cite[Lemma 1]{Bamb}}]
Let $\mathcal{O}$ be an $m$-ovoid of $\mathcal{P}_{r,e}$, then for any point $P \in \mathcal P_{r,e}$
$$
|P^{\perp} \cap \mathcal{O}|= \begin{cases}  (m - 1)(q^{r+e-2} + 1)+1 & \mbox{if } P \in \mathcal{O}, \\ m(q^{r+e-2} + 1) & \mbox{if } P \not \in \mathcal{O}.\end{cases}
$$
\end{lem}

An $m$-ovoid of $\mathcal{P}_{r,e}$ is said to be {\em trivial} if it is the empty set or consists of all points
of $\mathcal{P}_{r,e}$, and hence if $m=0$ or $m=\frac{q^r-1}{q-1}$, respectively. 
A 1-ovoid will be simply called {\em ovoid}. Let $\mathcal{O}$ and 
$\mathcal{O}^{\prime}$ be an $m$-ovoid and an $m^{\prime}$-ovoid of 
$\mathcal{P}_{r,e}$. If $\mathcal{O}$ and $\mathcal{O}^{\prime}$ are disjoint, 
then $\mathcal{O} \cup \mathcal{O}^{\prime}$ is an $(m+m^{\prime})$-ovoid of $\mathcal{P}_{r,e}$. Furthermore
the complement of an $m$-ovoid is a $\big(\frac{q^r-1}{q-1}-m\big)$-ovoid.

Ovoids, and $m$-ovoids for $m > 1$, are rare objects, see e.g.\ \cite{DeBeule,HirsThas} for an overview on 
known examples and some non-existence results for ovoids. Bounds on $m$ for the existence of $m$-ovoids of 
three classes of polar spaces can be found in \cite{Bamb}, with some recent improvements in 
\cite{DBMSxx}. 

\begin{defn}
Let $\mathcal{P}_{r,e}$ be a finite polar space of rank $r \geq 2$. A \emph{spread} 
$\mathcal{S}$ of $\mathcal{P}_{r,e}$ is a set of generators $\mathcal{P}_{r,e}$, 
which constitutes a partition of the pointset of $\mathcal{P}_{r,e}$.
\end{defn}

Spreads of polar spaces are a natural generalization of spreads
of projective spaces which received a lot of attention due to their importance for 
the theory of translation planes. Some polar spaces of rank 2 are dually isomorphic
to each other, hence results on spreads may be translated directly to results
on ovoids and vice versa. 

An $l$-system is a generalization of ovoids and spreads, introduced in \cite{ST1994} as well,
unifying the two concepts in one definition. 

\begin{defn}\label{lsy} Let $\mathcal{P}_{r,e}$ be a finite polar space of rank $r \geq 2$. 
An {\em $l$-system} of $\mathcal{P}_{r,e}$, with $0 \leq l \leq r-1$, is any set 
$\mathcal{S}=\{ \pi_1, \ldots, \pi_k \}$ of $k=q^{r+e-1}+1$ $l$-spaces of $\mathcal{P}_{r,e}$ 
such that no generator of $\mathcal{P}_{r,e}$ containing $\pi_i$ has a point in common with 
$\cup_{j \ne i} \pi_j$, $i=1,\ldots,k$.
\end{defn}  

By Definition \ref{lsy}, we have that a 0-system is an ovoid and a $(r-1)$-system is a spread. 
The following theorem shows that a putative $l$-system will yield an $m$-ovoid.

\begin{thm}[{\cite[6.1]{ST1994}}]
Let $\mathcal{S}$ be an $l$-system of the polar space $\mathcal{P}_{r,e}$, and $\mathcal{O}$ the set of points
covered by the $l$-dimensional subspaces in $\mathcal{S}$. Then $\mathcal{O}$ is an $m$-ovoid of $\mathcal{P}_{r,e}$
with $m = \frac{q^{l+1}-1}{q-1}$.
\end{thm}

There are known examples of $1$-systems of $Q^+(7,q)$, see e.g.\ \cite{LT2004,LT2005}, hence $(q+1)$-ovoids of 
$Q^+(7,q)$ exist as well. 

\subsection{Spreads of \texorpdfstring{$\mathrm{PG}(3,q)$ and ovoids of $Q^+(5,q)$}{PG(3,{\it q}) and ovoids of Q+(5,{\it q})}}\label{subsec: Klein}

Let $u = (u_1, u_2, u_3, u_4), v = (v_1, v_2, v_3, v_4)$ be two linearly independent 
vectors of $V = \mathbb{F}_q^4$. The {\em Plücker coordinates} of $(u, v)$ are $p_{ij} = u_iv_j - u_jv_i$ for every
$i \ne j$. Let $\tau$ be the map from pairs of linearly independent vectors of $V$ to points
of $\mathrm{PG}(5, q)$ defined by
$$
\tau(u,v)=\langle (p_{12}, p_{42}, p_{14}, p_{23}, p_{13}, p_{34}) \rangle.
$$
The map $\tau$ is well defined from the lines of $\mathrm{PG}(3,q)$ into the points of 
the quadric $Q^+(5,q):X_0X_5+X_1X_4+X_2X_3=0$ and it is called the {\em Klein correspondence}.
It follows that there is a bijective correspondence between the line spreads of 
$\mathrm{PG}(3,q)$\footnote{A {\em line spread} $\mathcal{S}$ of $\PG(3,q)$ is a set of $q^2+1$ 
lines such that each point of $\PG(3,q)$ belongs to exactly one line of $\mathcal{S}$.} and the 
ovoids of $Q^+(5,q)$. For more details on the Klein correspondence see \cite[\S 4.6]{Ball} or \cite[\S 4.8]{Beu}.

We can always assume that an ovoid of $Q^+(5, q)$ contains the points $(1,0,0,0,0,0)$ and 
$(0,0,0,0,0,1)$ and hence it can be written in the following form:
$$
\mathcal{O}(f_1,f_2)=\{(1,x,y,f_1(x,y),f_2(x,y),-yf_1(x,y)-xf_2(x,y)) \}_{x,y \in \mathbb{F}_q} \cup \{(0,0,0,0,0,1) \}
$$
for some functions $f_i:\mathbb{F}_q^2 \longrightarrow \mathbb{F}_q$ with $f_i(0,0)=0$, $i \in \{ 1,2 \}$; see \cite[p.\ 38]{Williams}.

\noindent The set $\mathcal{O}(f_1,f_2)$ is an ovoid if and only if 
\begin{eqnarray}\label{Eq:Intro2} \theta \big( (1,x_1,y_1,f_1(x_1,y_1),f_2(x_1,y_1),-x_1f_2(x_1,y_1)-y_1f_1(x_1,y_1)),\\ (1,x_2,y_2,f_1(x_2,y_2),
f_2(x_2,y_2),-x_2f_2(x_2,y_2)-y_2f_1(x_2,y_2)) \big) \ne 0,\nonumber
\end{eqnarray}
 
\noindent for every $(x_1,y_1) \ne (x_2,y_2)$ in $\mathbb{F}_q^2$, where $\theta$ is 
the reflexive form associated to the quadratic form of the quadric $Q^+(5,q)$. 

\begin{remark}
    Notice that in Equation (\ref{Eq:Intro2}), we did not consider the point $(0,0,0,0,0,1)$. This is not a restriction since it always holds that
   $$ \langle (1,x_1,y_1,f_1(x_1,y_1),f_2(x_1,y_1),-x_1f_2(x_1,y_1)-y_1f_1(x_1,y_1), (0,0,0,0,0,1)\rangle =1.$$
\end{remark}

\noindent Condition \eqref{Eq:Intro2} simplifies to
\begin{equation}\label{Eq:intro}
(x_1-x_2)(f_2(x_2,y_2)-f_2(x_1,y_1)) + (y_1-y_2) (f_1(x_2,y_2)-f_1(x_1,y_1)) \ne 0,\end{equation}
for every   $(x_1,y_1) \ne (x_2,y_2)$ in $\mathbb{F}_q^2$. 

Using the Klein correspondence, we can translate an ovoid $\mathcal{O}(f_1,f_2)$ 
of $Q^+(5,q)$ to the corresponding line spread of $\mathrm{PG}(3,q)$. This results in 
the set $\mathcal{S}(f_1,f_2)=\{ \ell_{x,y}: (x,y) \in \mathbb{F}_q^2 \} \cup \{ \ell_{\infty} \}$, 
where $\ell_{\infty}:X_0=X_1=0$ and
$$
\ell_{x,y}:\begin{cases} X_2=-f_1(x,y)X_0+f_2(x,y)X_1 \\ X_3=xX_0+yX_1. \end{cases}
$$

\section{Intersection patterns for \texorpdfstring{$m$}{\it m}-ovoids}\label{Patt}

In this section, we will assume that $n \geq 2$. 

\begin{lem}\label{le:ellipticToParabolic}
Suppose that $\mathcal{O}$ is an $m$-ovoid of $Q^-(2n+1,q)$, then it induces an $m$-ovoid in any $Q(2n,q) \subset Q^-(2n+1,q)$.
\end{lem}

\proof
Let $\pi$ be a hyperplane of $\mathrm{PG}(2n+1,q)$ such that $\pi \cap Q^-(2n+1,q)=Q(2n,q)$. Since any generator of $Q(2n,q)$ is also a generator of $Q^-(2n+1,q)$, it is easy to see that $\mathcal{O}' \coloneqq \pi \cap \mathcal{O}$ is an $m$-ovoid of $Q(2n,q)$.
\endproof

\begin{prop}\label{pattern}
Suppose that $\mathcal{O}$ is an $m$-ovoid of $Q^-(2n+1,q)$, with $m \geq 1$. Then any 
elliptic quadric $Q^-(2n-1,q) \subset Q^-(2n+1,q)$ meets $\mathcal{O}$ in exactly 
$(m-2)q^{n-1}+m$, $(m-1)q^{n-1}+m$ or $mq^{n-1}+m$ points, and all of these cases occur.
\end{prop} 

\proof
Consider an elliptic quadric $Q^-(2n-1,q)$ in an $(2n-1)$-dimensional space $\pi$, 
then $\pi^{\perp}$ is a line $\ell$ of $\mathrm{PG}(2n+1,q)$, meeting $Q^-(2n+1,q)$ in 2 points. 
Suppose that $|\ell \cap \mathcal{O}|=c$, then $c \leq 2$. Consider now any point $P \in \ell$. 
If $P \not \in Q^-(2n+1,q)$, then $P^{\perp} \cap Q^-(2n+1,q)=Q(2n,q)$ and hence $|P^{\perp} \cap \mathcal{O}|=m(q^n+1)$, 
since, by Lemma \ref{le:ellipticToParabolic}, $\mathcal{O}$ induces an $m$-ovoid in 
any $Q(2n,q) \subset Q^-(2n+1,q)$. If $P \in Q^-(2n+1,q) \setminus \mathcal{O}$, 
then $|P^{\perp} \cap \mathcal{O}|=m(q^n+1)$. For any point $P \in \mathcal{O}$, we 
get $|P^{\perp} \cap \mathcal{O}|=(m-1)q^n+m$.

By counting the pairs $\{ (P,R) : P \in \ell, R \in \mathcal{O} \setminus \ell, P \in R^{\perp}\}$, we find that
\[
m(q+1-c)(q^n+1)+c(m-1)(q^n+1)=m(q^{n+1}+1)-c+xq,
\]
This simplifies to
\[
m(q+1)(q^n+1)-cm(q^n+1)+cm(q^n+1)-c(q^n+1)=m(q^{n+1}+1)-c+xq,
\]
which yields $xq=mq(q^{n-1}+1)-cq^n$, with $x$ the number of points of $\mathcal{O} \cap \pi$ 
and $c \in \{0,1,2\}$. Thus the assertion follows.
\endproof

\begin{prop}
Suppose that $\mathcal{O}$ is an $m$-ovoid of $W(2n+1,q)$, $m \geq 1$. Then any symplectic 
space $W(2n-1,q) \subset W(2n+1,q)$ meets $\mathcal{O}$ in exactly $(m-c)q^{n-1}+m$ points for some $c \in \{ 0,1, \ldots, q+1\}$.
\end{prop}

\proof
Consider a symplectic space $W(2n-1,q)$ in an $(2n-1)$-dimensional space $\pi$, then $\pi^{\perp}$ 
is a line $\ell$ of $\mathrm{PG}(2n+1,q)$ contained in $W(2n+1,q)$. Suppose that $|\ell \cap \mathcal{O}|=c$, then $c \leq q+1$. 

If $P \in W(2n+1,q) \setminus \mathcal{O}$, then $|P^{\perp} \cap \mathcal{O}|=m(q^n+1)$. 
For any point $P \in \mathcal{O}$, we get $|P^{\perp} \cap \mathcal{O}|=(m-1)q^n+m$.

Hence, we can count the pairs $\{ (P,R): P \in \ell, R \in \mathcal{O}\setminus \ell, P \in R^{\perp} \}$ in order to find that
\[
m(q+1-c)(q^n+1)+c(m-1)(q^n+1)=m(q^{n+1}+1)-c+xq.
\]
This results in $x=m(q^{n-1}+1)-cq^{n-1}$, where $x$ is the number of points of $\mathcal{O} \cap \pi$ and  $c \in \{ 0,1, \ldots, q+1\}$. Thus the assertion follows.
\endproof

\begin{lem}\label{integer}
Let $m,c$ be two integers such that $m \geq 1$ and $c \in \{ 0,1,\ldots,q+1\}$. Then $(m-c)q^{n-1}+m=0$ if and only if $(n,c,m)=(2,q+1,q)$.
\end{lem}

\proof
We have that $(m-c)q^{n-1}+m=0$ if and only if $m=cq^{n-1}/(q^{n-1}+1)$. Note that 
$q^{n-1}+1 \mid cq^{n-1}$ if and only if $q^{n-1}+1 \mid c$. This is only possible whenever $(n,c)=(2,q+1)$. 
\endproof

This immediately implies the following corollary.

\begin{cor}
The following assertions hold:
\begin{itemize}
\item an $m$-ovoid $\mathcal{O}$ of $Q^-(2n+1,q)$, $m \geq 1$, has non-empty intersection with any $Q^-(2n-1,q) \subset Q^-(2n+1,q)$;
\item if an $m$-ovoid $\mathcal{O}$ of $W(2n+1,q)$, $m \geq 1$, has empty intersection with some $W(2n-1,q) \subset W(2n+1,q)$ then $(n,m)=(2,q)$. 
\end{itemize}
\end{cor}

\begin{prop}
Let $\mathcal{O}$ be an $m$-ovoid of $Q^-(2n+1,q)$, $m \geq 1$. If $\mathcal{O}$ contains an 
elliptic quadric $Q^-(2n-1,q)$, then $\mathcal{O}=Q^-(2n-1,q)$. 
\end{prop}

\proof
By contradiction, assume that $\mathcal{O} \ne Q^-(2n-1,q)$ and define $\Tilde{\mathcal{O}}$ 
to be the complement of $\mathcal{O}$. Then $\tilde{\mathcal{O}}$ is a 
$\Big( \frac{q^n-1}{q-1}-m\Big)$-ovoid that has empty intersection with $Q^-(2n-1,q)$, a contradiction. 
\endproof

\begin{prop}
Let $\mathcal{O}$ be an $m$-ovoid of $W(2n+1,q)$, $m \geq 1$. If $\mathcal{O}$ contains a 
symplectic space $W(2n-1,q)$ and $(n,m) \ne (2,q^2+1)$, then $\mathcal{O}=W(2n+1,q)$.
\end{prop}

\proof
Let $\tilde{\mathcal{O}}$ be the complement of $\mathcal{O}$. By contradiction, assume that 
$\mathcal{O} \ne W(2n+1,q)$. Then $\tilde{\mathcal{O}}$ is a $\Big(\frac{q^{n+1}-1}{q-1}-m \Big)$-ovoid 
and has empty intersection with $W(2n-1,q)$, a contradiction.
\endproof

\section{Constructing \texorpdfstring{$(q+1)$-ovoids of $Q^+(7,q)$}{({\it q+1})-ovoids of Q+(7,{\it q})}}\label{sec:construction}

We recall the following theorem obtained by B. Segre \cite{Segre} in 1965.

\begin{thm}[Segre's Theorem] 
Let $\mathcal{O}$ be a non-trivial $m$-ovoid of $Q^-(5,q)$, then $q$ is odd and $m=(q+1)/2$.
\end{thm}

We recall that the Klein correspondence is a bijective map between lines of $\PG(3,q^2)$ and 
points of $Q^+(5,q^2)$, see Section \ref{subsec: Klein}. In particular, the lines of a 
hermitian surface $H(3,q^2)$ are mapped onto points of an elliptic quadric $Q^-(5,q)$ 
obtained by intersecting $Q^+(5,q^2)$ with a subgeometry isomorphic to $\PG(5,q)$. Hence 
to any $m$-ovoid $\mathcal{O}$ of $Q^-(5,q)$ there corresponds a set of lines $\mathcal{L}$ 
of $H(3,q^2)$ with the property that every point of $H(3,q^2)$ is incident with exactly $m$ lines 
of $\mathcal{L}$. From Segre's Theorem it then follows that $q$ odd and $m=(q+1)/2$, 
thus the lineset $\mathcal{L}$ is called a {\em hemisystem}. In \cite{Segre}, Segre 
constructed an example of such a hemisystem of $H(3,9)$. Later, in \cite{Coss}, Cossidente 
and Penttila constructed a family of hemisystems of $H(3,q^2)$, for $q$ odd, hence 
generalizing the construction of Segre. From this we obtain that there exist 
examples of $(q+1)/2$-ovoids in $Q^-(5,q)$ for $q$ odd.

From now on, we will assume that $q$ is odd. Let $Q^+(7,q)$ be the hyperbolic 
quadric of $\mathrm{PG}(7,q)$. We describe a construction of a $(q+1)$-ovoid 
of $Q^+(7,q)$. Before doing this, we prove the following proposition.

\begin{prop}\label{prop:projection}
Consider $Q^+(7,q)$, and suppose that $\pi_1$ and $\pi_2$ are $5$-dimensional 
subspaces in the ambient space $\pg(7,q)$ such that $\dim(\pi_1\cap \pi_2)$ is a 3-space intersecting $Q^+(7,q)$ in a $Q^-(3,q)$. 
Suppose both $\pi_1$ and  $\pi_2$ intersect $Q^+(7,q)$ in an elliptic quadric $Q^-(5,q)$, 
that will be called $\mathcal Q_1$ and $\mathcal Q_2$ respectively. Then there exists a 
collineation $\Phi \in \PGO^+(8,q)$ such that $\mathcal{Q}_1$ is mapped on $\mathcal{Q}_2$ and the set
$\mathcal{Q}_1 \cap \mathcal{Q}_2$ is fixed pointwise.  
\end{prop}

\proof
Suppose that a $3$-dimensional subspace $\pi$ intersects $Q^+(7,q)$ in an elliptic quadric.
Then $\pi^\perp$, which is also of dimension $3$, intersects $Q^+(7,q)$ in another elliptic quadric. Hence both 
quadrics are an elliptic quadric $Q^-(3,q)$. Since $\pi$ and $\pi^\perp$ are disjoint, by abuse of notation, 
$\mathbb F_q^8 = \pi \oplus \pi^\perp$. 

Now suppose that $\varphi \in \PGO^-(4,q)$, i.e. is a similarity of $\mathcal{Q}_1$ in the ambient projective space 
$\pi^\perp$. Then we can extend $\varphi$ to an element $\tilde{\varphi}$ of $\PGO^+(8,q)$ as follows
\[
\tilde{\varphi}(x + y) = x + \varphi(y), \qquad x \in \pi, y \in \pi^\perp\,.
\]
The notation is explained as follows. The similarity $\varphi$ is represented 
by a $4 \times 4$ matrix $M$, and $\tilde{\varphi}(x + y)$ is represented by the $8 \times 8$ matrix
$$
M'=\begin{pmatrix} M & O \\ O & I_4 \end{pmatrix},
$$
where $O$ and $I_4$ are the zero-matrix and identity-matrix, respectively.
As such, $M'$ will represent indeed
a similarity of $Q^+(7,q)$, i.e. an element of $\PGO^+(8,q)$.

Now suppose that $\sigma$ is a 5-space through $\pi$, intersecting $Q^+(7,q)$ in an elliptic quadric $Q^-(5,q)$.
The polarity corresponding to this $Q^-(5,q)$ sends $\pi$ to $\ell = \sigma \cap \pi^\perp$.
Since this polarity is elliptic and $\pi$ intersects $Q^-(5,q)$ in an elliptic quadric $Q^-(3,q)$, 
$\ell$ will be a hyperbolic line with respect to $Q^-(3,q)$.
In general the span of a line $\ell \subset \pi^\perp$ and $\pi$ is hyperbolic, 
elliptic, or singular respectively if $\ell$ is elliptic, hyperbolic or singular. 
Since $\PGO^-(4,q)$ acts transitively on the hyperbolic lines in $\pi^\perp$, we now know that the 
pointwise stabilizer of $\pi \cap Q^+(7,q)$ in $\PGO^+(8,q)$ works transitively on the elliptic 5-spaces through $\pi$.
\endproof

\begin{construction}\label{con:ovoidHyperbolic}
First, we note that any $m$-ovoid of an elliptic quadric $Q^-(5,q)$ contained in $Q^+(7,q)$ is also an $m$-ovoid of $Q^+(7,q)$. 
Indeed, it is easy to see that any generator of $Q^+(7,q)$ meets $Q^-(5,q)$ in a generator of $Q^-(5,q)$. 
Let $\sigma$ be a generator of $Q^+(7,q)$ and let $\pi$ be a 5-dimensional subspace such that $Q^+(7,q) \cap \pi = Q^-(5,q)$. 
Now, $$\sigma \cap Q^-(5,q)= \sigma \cap \pi \cap Q^+(7,q)=\sigma \cap \pi.$$ Since $Q^-(5,q)$ has rank 2, then $\dim (\sigma \cap \pi) = 1$.
Moreover, $m$-ovoids of $Q^-(5,q)$ and $Q^+(7,q)$ have the same size $m(q^3+1)$.

Let $\pi_1$ and $\pi_2$ be two distinct 5-dimensional spaces such that $\pi_1 \cap Q^+(7,q)=\mathcal{Q}_1$ and $\pi_2 \cap Q^+(7,q)=\mathcal{Q}_2$ are two elliptic quadrics $Q^-(5,q)$ and $\dim( \pi_1 \cap \pi_2)=3$. It is possible to choose $\pi_1$ and
$\pi_2$ in such as way that $\pi_1 \cap \pi_2 \cap Q^+(7,q)= \mathcal{Q}_1 \cap \mathcal{Q}_2$ is an elliptic quadric $Q^-(3,q)$.
By  Proposition \ref{prop:projection}, there exists a collineation $\Phi \in \PGO^+(8,q)$ that maps $\mathcal{Q}_1$ on $\mathcal{Q}_2$ and fixes 
$\mathcal{Q}_1 \cap \mathcal{Q}_2$ pointwise.

Now, let $\mathcal{O}_1$ be a $(\frac{q+1}{2})$-ovoid of $\mathcal{Q}_1$. Then, $\Phi(\mathcal{O}_1)$ is a $(q+1)/2$-ovoid of $\mathcal{Q}_2$, hence its complement $\mathcal{O}_2$ in $\mathcal{Q}_2$ is again a $(q+1)/2$-ovoid. Since $\mathcal{O}_1$ and $\mathcal{O}_2$ are disjoint in $Q^+(7,q)$, then $\mathcal{O}_1 \cup \mathcal{O}_2$ is a $(q+1)$-ovoid of $Q^+(7,q)$.
\end{construction}

 The following Theorem summarizes the result of this Section.
\begin{thm}
    There exist $(q+1)$-ovoids in $Q^+(7,q)$ which are the union of two $\frac{q+1}{2}$-ovoids in distinct $Q^-(5,q)$'s.
\end{thm}
\proof
See construction \ref{con:ovoidHyperbolic}.
\endproof

\begin{remark}
A natural question is whether this ``glueing'' approach could be feasible in polar spaces of other type than (hyperbolic) quadrics. 
Although it cannot be excluded that a similar approach could work for symplectic or hermitian spaces, 
it should be noted that polar spaces $Q^-(5,q) \subset Q(6,q) \subset Q^+(7,q)$ have the exceptional property that an $m$-ovoid of $Q^-(5,q)$
is also an $m$-ovoid of $Q(6,q)$ and also of $Q^+(7,q)$. This series of three embedded polar spaces, with the property that 
any generator of a larger one meets a smaller one also in a generator, only occurs in the family of quadrics. 
It might be worth investigating the generalization into arbitrary dimension, i.e. based on the series 
$Q^-(2n-1,q) \subset Q(2n,q) \subset Q^+(2n+1,q)$, $n > 3$. The main difficulty however will be to find suitable $m$-ovoids of $Q^-(2n-1,q)$,
$n > 3$.
\end{remark}
 
\section{Spreads without tangents to \texorpdfstring{$Q^-(3,q)$}{Q-(3,{\it q})}}\label{notang}

Consider the elliptic quadric $Q^-(3,q)$ with ambient projective space $\mathrm{PG}(3,q)$. 
Our first goal is to construct a line spread of $\mathrm{PG}(3,q)$ containing as 
many 2-secants\footnote{A line $\ell$ is said to be {\em $i$-secant} to 
a set $\Gamma$ if it has exactly $i$ points in common with $\Gamma$. 
In particular, if $i=1$ the line is also called {\em tangent}.}
to $Q^-(3,q)$ as possible. The maximum number clearly equals $|Q^-(3,q)|/2=(q^2+1)/2$, 
which is evidently only attainable if $q$ is odd.
Note that this number is reached if and only if the line spread contains no tangent lines to $Q^-(3,q)$.
We will show that this number 
can be reached if $q=3^h$, $h$ odd, using the Klein
correspondence (Section~\ref{subsec: Klein}), or when $q \equiv 1 \pmod 4$. 

\noindent \underline{Case 1: $q=3^h$, $h$ odd.}\\
Then $q \equiv -1 \pmod 4$ and so $-1$ is a non-square.

Let $f_1(x,y)=x+y$, $f_2(x,y)=x+2y$. Then, Equation \eqref{Eq:intro} reads
\[
(x_1-x_2)^2+(y_1-y_2)^2 \ne 0.
\]
Since $-1$ is not a square in $\mathbb{F}_q$, the equation above is satisfied 
for any $(x_1,y_1) \ne (x_2,y_2)$ in $\mathbb{F}_q^2$. Hence, we get the following.
\begin{thm}
Let $q=3^h$, $h$ odd, and suppose that $f_1(x,y)=x+y$, $f_2(x,y)=x+2y$. Then 
$\mathcal{O}(f_1,f_2)$ is an ovoid of $Q^+(5,q)$ and the set 
$\mathcal{S}(f_1,f_2)=\{ \ell_{x,y}: (x,y) \in \mathbb{F}_q^2 \} \cup \{ \ell_{\infty} \}$ 
is a line spread of $\mathrm{PG}(3,q)$, where $\ell_{\infty}:X_0=X_1=0$ and
$$
\ell_{x,y}:\begin{cases} X_2=-(x+y)X_0+(x+2y)X_1 \\ X_3=xX_0+yX_1. \end{cases}
$$
\end{thm}

The quadratic form $F(X_0,X_1,X_2,X_3)=X_0X_1+X_2^2+X_3^2$ defines an elliptic quadric.
Indeed, since $-1$ is not a square, $X_2^2 + X_3^2$ is irreducible.
Now, consider the line spread $\mathcal{S}(f_1,f_2)$ defined in the previous theorem. 
Note that the line $\ell_{\infty}$ is a 0-secant to $Q^-(3,q)$ and the line $\ell_{0,0}$ is 
a 2-secant to $Q^-(3,q)$. Let $(x,y)\ne (0,0)$, then the line $\ell_{x,y}$ contains the points 
$(1,0,-x-y,x)$ and $(0,1,x+2y,y)$, which do not lie on $Q^-(3,q)$. Hence, this line is tangent if and only if 
\begin{eqnarray}
F((1,0,-x-y,x)+T(0,1,x+2y,y))= \nonumber \\ (x^2+xy+2y^2)T^2+(x^2+2xy+2y^2+1)T+2x^2+2xy+y^2 \nonumber
\end{eqnarray}
has a unique root, which is equivalent to
$$
x^4+2x^2y^2+y^4+x^2+2xy+2y^2+2=0.
$$
So finding tangent lines to $Q^-(3,q)$ is equivalent to finding $\fq$-rational 
points of the curve $\mathcal{C}:H(x,y)=0$, where $H(x,y)=x^4+2x^2y^2+y^4+x^2+2xy+2y^2+2$. 

Let $i \in \overline{\mathbb{F}}_q$ be 
such that $i^2=-1$. Let $t\coloneqq x+iy$ and $z\coloneqq x-iy$. Now, the equation of $\mathcal{C}$ reads
$$
t^2z^2+(i+2)t^2+(2i+2)z^2+2=0,
$$
or equivalently
\begin{eqnarray}\label{fq3:1}
(z^2-(1-i))(t^2-(1+i))=0.
\end{eqnarray}
Note that the smallest field containing the square roots of $1-i$ and $1+i$ is $\mathbb F_{3^4}$.
Since $q = 3^h$, with $h$ odd, $\mathbb F_{3^4}$ is not a subfield of $\mathbb F_{q^2}$, thus Equation \eqref{fq3:1} has no solutions over $\mathbb F_{q^2}$.
It follows that the 
curve $\mathcal{C}$ has no $\mathbb{F}_{q}$-rational points.
Therefore, the spread $\mathcal{S}(f_1,f_2)$ of $\mathrm{PG}(3,q)$ has no tangent lines to $Q^-(3,q)$.

\noindent \underline{Case 2: $q \equiv 1 \pmod 4$.}\\
We construct a Desarguesian spread, and show that it does not contain any tangent 
lines to $Q^-(3,q)$. Choose a non-square $\alpha \in \fq$. Then the quadratic 
form $F(X_0,X_1,X_2,X_3)=X_0X_1+X_2^2-\alpha X_3^2$ defines an elliptic quadric.

Every line in $\mathrm{PG}(3,q)$ can be represented as the column space of a 
full rank $4 \times 2$ matrix. Write $A=\begin{pmatrix} 0 & \alpha \\ 1 & 0 \end{pmatrix}$, $I_2=\begin{pmatrix} 
 1 & 0 \\ 0 & 1 \end{pmatrix}$ and $O = \begin{pmatrix} 0 & 0 \\ 0 & 0 \end{pmatrix}$, then
\[
\mathcal{S}=\Bigg\{ \begin{pmatrix} I_2 \\ xI_2 + yA \end{pmatrix} : x,y \in \fq \Bigg\} \cup \Bigg\{ \begin{pmatrix} O \\ I_2 \end{pmatrix} \Bigg\}
\]
is a Desarguesian spread. Note that the line $\begin{pmatrix} O \\ I_2 \end{pmatrix}$ is 
0-secant and the line $\begin{pmatrix} I_2 \\ 0 \end{pmatrix}$ is 2-secant.

So take another line $\begin{pmatrix} I_2 \\ xI_2 + yA \end{pmatrix} \in \mathcal{S}$, $(x,y) \ne (0,0)$. 
This line goes through the points $(1,0,x,y)$ and $(0,1,\alpha y,x)$ which do not lie on $Q^-(3,q)$ due to $(x,y) \ne (0,0)$. 
Hence this line is tangent to $Q^-(3,q)$ if and only if
\[
F((1,0,x,y)+T(0,1,\alpha y,x))=\alpha(\alpha y^2-x^2)T^2+T+(x^2-\alpha y^2)
\]
has a unique root, which is equivalent to
\[
1-4 \alpha(\alpha y^2-x^2)(x^2-\alpha y^2)=1+4\alpha (x^2-\alpha y^2)^2=0,
\]
since $\alpha(\alpha y^2-x^2)\ne 0$ if $(x,y) \ne (0,0)$. Note that since $\alpha$ is a 
non-square, $4\alpha (x^2-\alpha y^2)^2$ is a non-square. If $q \equiv 1 \pmod 4$, $-1$ is a 
square and the above equation has no solutions. Therefore, the Desarguesian spread $\mathcal{S}$ 
contains no tangent lines to $Q^-(3,q)$.

The results above are summarized in the following proposition, where we use the fact that $3^h \equiv 1 \pmod 4$ if $h$ is even.

\begin{prop}\label{prop:maximalSpread}
Consider $Q^-(3,q)$ in the ambient space $\pg(3,q)$. If $q=3^h$ or $q\equiv 1 \pmod 4$,
then there exists a line spread of $\pg(3,q)$ containing $|Q^-(3,q)|/2=(q^2+1)/2$ $2$-secant lines to $Q^-(3,q)$.
\end{prop}

The existence of spreads with no tangents to $Q^-(3,q)$ will be used to construct $m$-ovoids in $Q^+(7,3)$ for $m \in \{2,4,6,8,10\}$. 
An extended computer search could lead to examples of $m$-ovoids of $Q^+(7,3)$ for $m \in \{2,4,6,8,10\}$ as well. However,
given an example of an $m$-ovoid found using a computer search, it would not be easy to show that it is the disjoint union of 
$m'$-ovoids for $m' < m$. Therefore a construction result as presented in Theorem~\ref{Thm:m-OvoidsForQ=3} for $m$-ovoids
is rather interesting. We need one more lemma.

\begin{lem}
 \label{Lem:2-ovoidTroughEveryPair}
 Consider the elliptic quadric $Q^-(5,3)$ and let $\sigma$ be an elliptic $3$-space.
 Then for any 2 distinct points $P,R \in \sigma$, there exists a 2-ovoid $\mathcal O$ 
 of $Q^-(5,3)$ with $\mathcal O \cap \sigma = \{P,R\}$.
\end{lem}

\begin{proof}
There exists a 2-ovoid $\mathcal O_1$ of $Q^-(5,3)$.
Let $\pi$ be the ambient 5-space of $Q^-(5,3)$.
By \Cref{pattern} there is a $3$-space $\sigma' \subset \pi$ which intersects $\mathcal O_1$ in $2$ points.
Since $\PGO^-(6,3)$ acts transitively on the elliptic 3-spaces, there also exists 
a 2-ovoid $\mathcal O_2$ intersecting $\sigma$ in 2 points $\{P',R'\}$.
The group $\PGO^-(4,3)$ acts 2-transitively on the points of $Q^-(3,q)$.
Thus there exists a collineation $\varphi$ of $\sigma$ which fixes $Q^-(5,3) \cap \sigma$ 
and maps $\{P',R'\}$ to $\{P,R\}$. As in \Cref{prop:projection}, we can choose an 
appropriate element of the stabilizer of $\sigma^\perp \cap Q^-(5,3)$, which is 
isomorphic to $\PGO^+(2,3)$, to extend $\varphi$ to an element of $\PGO^-(6,3)$.
This extended collineation maps $\mathcal O_2$ to a 2-ovoid of $Q^-(5,3)$ 
intersecting $\sigma$ exactly in $\{P,R\}$.
\end{proof}

\begin{thm}\label{Thm:m-OvoidsForQ=3}
The hyperbolic quadric $Q^+(7,3)$ contains $m$-ovoids, for $m \in \{2,4,6,8,10\}$, consisting of the disjoint union of $2$-ovoids.
\end{thm}

\begin{proof}
We will prove that $Q^+(7,3)$ has 5 pairwise disjoint 2-ovoids.
The theorem then immediately follows by taking unions of some of these 2-ovoids.

Take an elliptic $3$-space $\sigma$.
Then $\sigma^\perp$ is also an elliptic 3-space.
By \Cref{prop:maximalSpread}, there exists a line spread of $\sigma^\perp$ 
containing 5 2-secant lines $\ell_1, \dots, \ell_5$ to $Q^+(7,3) \cap \sigma^\perp \cong Q^-(3,3)$.
Let $\pi_i$ denote the span of $\ell_i$ and $\sigma$ for $i = 1, \dots, 5$.
Then $\pi_1, \dots, \pi_5$ are elliptic 5-spaces pairwise intersecting in $\sigma$.
Partition the 10 points of $\sigma \cap Q^+(7,3)$ into 5 pairs $\{P_i,R_i\}$, $i = 1, \dots, 5$.
For every $i$ choose a 2-ovoid $\mathcal O_i$ of $\pi_i$ which intersects $\sigma$ precisely in $\{P_i,R_i\}$.
Such a 2-ovoid exists by \Cref{Lem:2-ovoidTroughEveryPair}.
Then $\mathcal O_1, \dots, \mathcal O_5$ are pairwise disjoint 2-ovoids of $Q^+(7,3)$.
\end{proof}

\begin{remark}
Notice that the same technique cannot be used for $q>5$, since, by \Cref{pattern}, 
every $\frac{q+1}2$-ovoid of $Q^-(5,q)$ intersects every elliptic 3-space in more than $|Q^-(3,q)|/3$ points.
Finally, using a non-exhaustive computer search, we found no indication of the existence of $3$ pairwise disjoint 
$3$-ovoids in $Q^+(7,5)$.
\end{remark}

\section*{Acknowledgement}
The research of Giovanni Giuseppe Grimaldi was supported by
the Italian National Group for Algebraic and Geometric Structures and their
Applications (GNSAGA - INdAM).
The authors also want to thank the referees for their valuable suggestions to improve this article.

\end{document}